
\input amstex
\documentstyle{amsppt}
\magnification=\magstep1
\toc \nofrills{Table of Contents}
\widestnumber\head{10.}
\head 1.  Introduction$\dotfill 1$\endhead
\head 2.  The Structure Theorem$\dotfill 3 $\endhead
\head 3.  A Corollary of the Structure Theorem$\dotfill 8 $\endhead
\head 4.  $(r_m, \delta_m)$ as a random walk$\dotfill 9 $\endhead
\head 5.  The ensemble $\Phi_m^{-1} (r_m, \delta_m)$$\dotfill 11 $\endhead
\vskip .3cm
\specialhead {} References $\dotfill 15 $\endspecialhead
\endtoc
\vglue -.01pc
\topmatter
\leftheadtext{Ya. G. Sinai}
\rightheadtext{Statistical $(3x +1)$ -- problem}
\TopOrBottomTagsOnSplits
\TagsOnRight
\nologo
\pageheight{7.0in}
\pagewidth{5.1in}
\title {Statistical $\bold{(3x+1})$ -- problem }\endtitle
\author{Ya. G. Sinai$^\star$}\endauthor
\footnotetext"{$^\star$}"{Mathematics Department of Princeton University and Landau Institute of Theoretical Physics}
\dedicatory{Dedicated to the memory of J. Moser}\enddedicatory
\parskip=5.5pt
\NoBlackBoxes
\endtopmatter

\define\jK1{\Delta_{\tsize{k^{(0)},\ldots ,k^{(j+1)}}}^{(j)}}

\define\Dd2{\Delta_{\tsize{k^{(0)},\ldots ,k^{(i-2)}}}^{(i-1)}}
\define\di+{\Delta_{\tsize{k^{(0)},\ldots ,k^{(i-2)}}}^{(i-1)}}

\define\rjdj{r_j , \delta_j}
\define\kme{k_1, \ldots , k_m , \epsilon}
\define\1m{m+1}
\define\2m{m+2}

\define\({\left(}
\define\){\right)}

\document

\baselineskip=13pt

\head{\bf \S1. Introduction}\endhead

Take an odd number $x >0$.  Then $3x+1$ is even and one can find an integer
$k> 0$ so that $y= {3x+1\over 2^k}$ is again odd.  We get in this way the
mapping $T, T x=y$.  It is clear that except being odd $y$ is also not
divisible by 3.  By this reason the natural domain for $T$ is the
set $\sqcap$ of positive $x$ not divisible by 2 and 3.
The point $x =1$ is the fixed point of $T$ and it is the
famous $(3x+1)$--problem which asks whether it is true that for
every $x \in \sqcap$ one can find $n(x)$ such that $T^{n(x)} x=1$.
The best references concerning this problem are the expository
paper by J. Lagarias \cite{L} and the book by G. Wirsching \cite{W},
see also the annotated bibliography on $(3x +1)$--problem
prepared by J. Lagarias.
There one can find a lot of information about the history of
the problem and its various modifications.
We call the statistical $(3x+1)$--problem the basic question for 
$x$ belonging to a subset of density 1.  In this paper we discuss
some version of the statistical $(3x+1)$--problem.

\pagebreak

The main result of this paper is the following Structure Theorem
which we formulate now and prove in \S2.

We have $\sqcap = 1\cup \sqcap^{+1} \cup \sqcap ^{-1}$ where
$\sqcap^{+1} = \{ 6p+1 , p > 0\}$,
$\sqcap^{-1} = \{ 6p-1, p>0\}$.
Let us write if in the definition of $T$
we divide $3x+1$ 
by $2^k$.  Fix $k_1 , k_2 , \ldots , k_m$,
$k_j >0$ and integer, $1 \leq j \leq m$.  We ask what is the set
of $x \in \sqcap$ to which one can successively apply
$T^{(k_1)} , T^{(k_2)} , \ldots , T^{(k_m)}$.
The Structure Theorem gives the answer.

\proclaim{Structure Theorem}
Let $k_1,\ldots , k_m$ be given and $\epsilon_m = \pm 1$.
The set of $x \in \sqcap^{\epsilon_m}$ to which one can apply
$T^{(k_1)}, T^{(k_2)} , \ldots , T^{(k_m)}$ is an arithmetic
progression $\Sigma^{(k_1,\ldots , k_m, \epsilon_m)}$
$= \{6 \cdot (2^{k_1+k_2+\cdots+ k_m} p + q_m) + \epsilon_m \}$ for some
$q_m = q_m(k_1 , \ldots , k_m , \epsilon_m)$,
$1 \leq q_m \leq 2^{k_1 + \cdots + k_m}$.  The image 
$T^{(k_m)} \cdot T^{(k_{m-1})}  \cdots T^{(k_1)}$
$(\Sigma^{(k_1 , \ldots , k_m , \epsilon_m)}) = \wedge_{r_m , \delta_m}^{(m)}$
$=\{6(3^m p+r_m)+\delta_m\}$ for some 
$r_m = r_m (k_1, \ldots ,  k_m, \epsilon_m)$,
$1 \leq r_m \leq 3^m$, and $\delta_m = \delta_m (k_1 , \ldots , k_m, \epsilon_m ) = \pm 1$.  Even more, for each $p > 0$
$T^{(k_m)} , T^{(k_{m-1})}  \cdots  T^{(k_1)}$
$(6(2^{k_1 + \cdots + k_m} p +q_m ) + \epsilon_m) = 6(3^mp+r_m)+ \delta_m$
with the same $p$.
\endproclaim

The proof of this theorem goes by induction.  First we check
the statement for $m=1$ and then derive it for $m+1$ assuming that
it is true for $m$.
A. Kontorovich has a shorter proof of this theorem.

This theorem plays the role of symbolic representation in dynamics.

In \S3 we prove a simple statistical statement which follows directly
from the Structure Theorem.  Take $x_0 \in \sqcap$,
$x_m = T^m x_0$, $y_m = \ell n \; x_m$ and
$z_m = y_m - y_0$.  Assume that $1 \leq m \leq M$ and
$$
\omega \left({m\over M}\right) =
   {z_m + m (2 \ell n \; 2 - \ell n \; 3)\over \sqrt {M}} .
$$
We show that $\omega(t), 0 \leq t \leq 1$, behave as
Wiener trajectories.  More precisely, let $M= 2^n$,
$n \rightarrow \infty , \; \tau_1 = {t_1\over 2^{n}}$,
$\tau_2 = {t_2\over 2^{n}} , \ldots , \tau_s = {t_s\over 2^{n}}$
where $t_1 , t_2 , \ldots , t_s$ are integers,
$0 \leq t_j < 2^{n}$, $1 \leq j \leq s$,
and $\tau_1 , \ldots , \tau_s$ are fixed.
Consider the following probability
$$
P_n = P \{ x_0 | a_1 \leq \omega (\tau_1) \leq  b_1 , \ldots ,
a_s \leq \omega (\tau_s) \leq b_s \}
$$
Here, $a_1, b_1 , \ldots , a_s , b_s$ are fixed numbers.
The probability $P$ of a set is understood as the density 
of this set wrt $??$ provided that the density exists. 
The following theorem holds.

\proclaim{Theorem 3.1}
The probability $P_n$ tends as $n \rightarrow \infty$ to the
probability given by the Wiener measure with the zeroth drift and some
diffusion constant $\sigma > 0$. 
\endproclaim

In \S3 we prove it for $s = 1$.  General case can be
obtained in a similar way.

An analogous theorem was proven
recently by K. A. Borovkov and
D. Pfeifer (see \cite{BP}).

In \S4 we study some properties of $r_m (k_1 , \ldots , k_m , \epsilon_m)$.
Technically the most important part is in \S5 where we analyze the
ensemble of those $(k_1, k_2, \ldots , k_m , \epsilon_m)$
for which $k_1 + k_2 + \cdots + k_m = k$ and
$r_m (k_1, k_2, \ldots , k_m, \epsilon_m)$
$= r_m , \delta_m (k_1, \ldots , k_m , \epsilon_m) = \delta_m$ are fixed.


During the work on $(3x+1)$--problem I had many discussions
with A. Bufetov,\newline
\noindent  N. Katz,  A. Kontorovich, L. Koralov, J. Lagarias,
and T. Suidan.  It is my pleasure to thank all of them for very
useful contacts. 
The criticism from J. Lagarias was especially important.  He found
a serious gap in the previous version of this paper.
I am planning to discuss the related questions in a
forthcoming publication.
I thank also a referee for many useful remarks.
The financial supports from NSF, grant DMS--0070698
and RFFI, grant N99--01--00314 are highly appreciated.

\head{\bf \S2.  The Structure Theorem}\endhead

The formulation of the Structure Theorem was given
in \S1.  Here we give its proof.  First we consider the
case $m=1$.

Assume that $\epsilon_m = +1$ and take $x=6p+1 \in \sqcap^{+1}$.
For given $k_1$ we should have
$$
3x+1 = 18p +4 = 2^{k_1} (6t + \delta_1)
$$
for some $\delta_1 = \pm 1$.
\roster
\item"{$a_1$)}"  $k_1 = 1$.  Then
$$ 9p+2 = 6t + \delta_1.$$ 
\endroster
This shows that $p $ has to be odd, $p=2p_1+1$ and
$$
6(3p_1 +2) -1 = 6t+ \delta_1 .
$$
Therefore in this case $\delta_1  = -1, t= 3p_1 +2$,
i.e. $q_1 (1, +1) = 1 , \; r(1, +1) = 2$.
\roster
\item"{$a_2$)}"  $k_1 > 1$.  In this case
$$ 9p+2 = 2^{k_1 -1} (6t+\delta_1) .$$
\endroster
This shows that $p$ has to be even, $p = 2p_1$ and
$$
6 \cdot 3 p_1= 6\cdot 2^{k_{1}-1}
   t+\delta_1 2^{k_1 -1 }-2 \,.
\tag"{(2.1)}"
$$
The number $\delta_1 2^{k_1 -1} -2$ is even. 
The value $\delta_1$ should be chosen so that 
$\delta_1 \cdot 2^{k_1 -1} -2 \equiv 0$ $\pmod 3$.
If $k = 3, 5, 7, \ldots $ then $\delta_1$ should be $-1$.
If $k = 2, 4, 6, \ldots $ then $\delta_1 = 1$.  In other words,
$\delta _1 2^{k_1} \equiv 1 \pmod 3$.
Since $\delta_1 = \delta_1^{-1}$ the last expression takes the form
$$
\delta_1 \equiv 2^{k_1} \pmod 3 \, ,
\tag"{(2.2)}"
$$
i.e. $\delta_1$ is uniquely determined by $k_1$.

Returning back to (2.1) we have
$$
3p_1 - 2^{k_1 -1} t = {\delta_1 2^{k_1 -2} -1\over 3} \,.
$$

Let us write $p_1 = 2^{k_1 -1} s+ \bar{q}_1 , t = 3s+ \bar{r}_1$.
For $\bar{q}_1 , \bar{r}_1$ we have the equation
$$
3 \bar{q}_1 - 2^{k_1 -1} \bar{r}_1 =
   {\delta_1 2^{k_1 -2} -1\over 3} \,.  
\tag"{(2.3)}"
$$
If $\delta_1 =1$ then $rhs$ of (2.3) is non--negative
and less than $2^{k_1 -1}$.
If $k_1 =2$
 then it is zero and the solution to
the equation (2.3) takes the form
$\bar{q}_1 =2, \bar{r}_1 = 3$, i.e. $ q(2,+1) = 4$,
$r_1 (2, +1) =3$.

If $k_1 > 2$ then $rhs$ of (2.3) is positive.
Consider the abelian group $Z_2^{k_{1 -1}}$ of numbers $\mod 2^{k_1-1}$.
The multiplication by 3 is an automorphism of this group.

(2.3) is the equation for $\bar{q}_1$ in $Z_2^{k_1 -1}$ and
it has a unique solution.  The value $\bar{r}_1 , 1 \leq \bar{r}_1 \leq 3$
is also defined uniquely.  Thus $q (k_1, +1) = 2\bar{q}_1$,
$r(k_1 , +1) = \bar{r}_1$.

If $\delta_1 = -1$ then $rhs$ of (2.3) is negative.
We rewrite (2.3) as follows:
$$
3 \bar{q}_1 -2^{k_1-1} (\bar{r}_1 -1) = 2^{k_1-1} +
   {\delta_1 2^{k_1-2} -1\over 3}\,.
\tag"{(2.4)}"
$$
Now $rhs$ is positive and we can use the same arguments
as before to find $0 < \bar{q}_1 \leq 2^{k_1-1} $, 
$0 \leq \bar{r}_1 -1 < 3$.
Therefore $r({k_1} , +1) = \bar{r}_1,  1 \leq \bar{r}_1 \leq 3$.

The case $x \in \sqcap^{-1}$ is considered in a similar way.
We write $x = 6p-1, p > 0$.  For given $k_1$
we have the equation
$$
3x+1 = 18p -2 =18 p^\prime +16 = 2^{k_1} (6t + \delta_1 )
$$
for $p^\prime = p-1 \geq 0$ and some $\delta_1 = \pm 1$.

$a_1$)  $k_1 = 1$.  Then
$$
9 p^\prime +8 = 6t + \delta_1 \,.
$$
This shows that $p^\prime$ has to be odd, $p^\prime = 2s+1$ and
$$
18s - 6t = - 17 + \delta_1 \, .
$$
Therefore $\delta_1 = -1$ and $3s +3 = t$.
From the last expression $p^\prime = 2s +1$,
i.e. $q(2, -1) =1$, $r(2,-1) = 3$ and $s=0$.

$a_2)$  $k_1 > 1$.  Then
$$
9p -2 = 2^{k_1 -1} ( 6 t + \delta_1) \,.
$$
$p$ has to be even, $p = 2p_1, p_1 > 0$ and
$$
9p_1 -2^{k_1 -2} \cdot 6t = 2^{k_1-1} \delta_1 -8 \,.
$$
The last expression shows that $p_1$ also has to be
even, $p_1 = 2p_2 , p_2 > 0$ and
$$
9 p_2-3 \cdot 2 ^{k_1-2} t = 2^{k_1-2} \delta_1 -4 \,.
$$
$Rhs$ must be divisible by 3.  This gives 
$2^{k_1 -2} \delta_1 -4 \equiv 0 \pmod 3$ or
$2^{k_1} \equiv \delta_1 \pmod 3$.
We get the equation
$$
3p_2 - 2^{k_1-2} t = 
  {2^{k_1-2} \delta_1 -4\over 3}
\tag"{(2.5)}"
$$
If $k_1 = 4, \delta_1 = 1$, then $rhs$ of (2.5) is zero
and $p_2 = 2^{k_1 -2} s, t= 3s, s> 0$.
In other words, $p = 2^4s, t = 3s, q_1 (4, -1) =0$,
$r_1 (4, -1) = 0$.
In order to comply with the formulation of the theorem we
change our choice to $q_1(4,-1) = 2^4$, $r_1(4,-1) = 3$ and
$p \geq 0$.

If $k_1 = 2 $ then $\delta_1 =1$ and it is easy to check that
$q(2,-1)=3, r(2,-1) =2$.

If $k_1 > 4, \delta_1 =1$ we argue as before.
$Rhs$ of (2.5) is positive and less than $2^{k_1 -2}$.
We put $p_2 = 2^{k_1 -2} s + \bar{p}_2$, $t= 3s + \bar{t}$
and for $\bar{p}_2, \bar{t}$ we get the equation
$$
3 \bar{p}_2 - 2^{k_1-2} \bar{t} =
  {2^{k_1-2} \delta_1-4\over 3}
\tag"{(2.6)}"
$$
which has the unique solution $\bar{p}_2 , 0 < \bar{p}_2 \leq 2^{k_1-2}$,
and $\bar{t}$, $1 \leq \bar{t} \leq 3$. 
This gives $p = 2^{k_1 s} +4 \bar{p}_2$,
$t= 3s+\bar{t}$, i.e. $q_1(k_1 , -1) = 4 \bar{p}_2$,
$r_1(k_1 , -1) = \bar{t}$.

If $\delta_1 = -1$ and $rhs$ of (2.6) is negative we modify it as before
$$
3\bar{p}_2 -2^{k_1 -2} (\bar{t}-1) = 2^{k_1-2}
   + {2^{k_1-2} \delta_1 -4\over 3}
\tag"{$(2.6^\prime)$}"
$$
Now $rhs$ of (2.6$^\prime$) is positive and we can find
$1 \leq \bar{p}_2 \leq 2^{k_1-2}$, $0 \leq \bar{t} -1 < 3$
satisfying (2.6$^\prime$).  Then $p = 2^{k_1} s+4 p_2$,
$q(k_1 , -1) = 4\bar{p}_2$ and
 $r_1(k_1, -1) = \bar{t}, 1 \leq \bar{t} \leq 3$.

The case $m > 1$ is considered by induction.
Suppose that for some $m \geq 1$ the Structure Theorem
is proven, i.e.
$T^{(k_m)} \cdot T^{(k_{m-1})} \cdots T^{(k_1)} (6(2^{k_{1}+\cdots+ k_m} s)$
$+q_m(k_1, \ldots, k_m, \epsilon_m)) = 6(3^m s + r_m ( k_1,\ldots , k_m , \epsilon_m))+\delta_m (k_1, \ldots , k_m , \epsilon_m)$.
Denote $x= 6(3^m s + r_m) + \delta_m$.
Then $3x+1 = 2^{k_{m+1}} y$ where $y$ is odd and
$$
6(3^{m+1} s +3r_m ) + 3\delta_m +1 =2^{k_{m+1}} y
$$
or
$$
3 \cdot 3^{m+1} s + 9 r_m +
  {3\delta_m +1\over 2} = 2^{k_{m+1}-1} y \,.
\tag"{(2.7)}"
$$
If $k_{m+1} = 1, \delta_m =1$ 
then (2.7) takes the form
$$
3 \cdot 3^{m+1} s + 9 r_m + 2 = y \, .
$$

If $r_m$ is even, $r_m = 2 r_m^{(1)}$, then $s$ must be odd,
$s = 2s_1 +1$,
$$
3 \cdot 3^{m+1} (2s_1 , +1) + 9\cdot 2 r_m^{(1)} + 2 = y \,,
$$
or
$$
6\left(3^{m+1} s_1 + 3 r_m^{(1)} + {3^{m+1} +1\over 2} \right) -1 = y \,.
$$
This shows that 
$\delta_{m+1} = -1, r_{m+1}$ 
$(k_1, \ldots , k_{m+1}, \epsilon_{m+1}) = 3r_m^{(1)} + {3^{m+1} +1\over 2}$,
$q_{m+1} = g_m +2^{k_1 +\cdots + k_m}$.

If $r_m$ is odd, $r_m = 2 r_m^{(1)} +1$ then $s$ has to be even,
$s = 2s_1$ and
$$
6(3^{m+1} s_1 + 3 r_m^{(1)} +2) - 1 = y \,.
$$

We conclude that $\delta_{m+1} =-1$,
$r_{m+1} (k_{1, \dots ,} k_{m+1}, \epsilon_{m+1} )$
$= 3 r_m^{(1)} +2, q_{m+1} = q_m$.

The case $k_{m+1} =1, \delta_m = -1$ is considered in a similar way.

Now let $k_{m+1} > 1$.  If $r_m$ is even,
$r_m = 2r_m^{(1)}$ and $\delta_m =1$ then $s$ has to be even,
$s = 2s_1$ (see (2.7)) and from (2.7)
$$
6\cdot 3^{m+1} s_1 + 18 r_m^{(1)} +2 = 2^{k_{m+1} -1} y \,.
\tag"{(2.8)}"
$$
If $y=6t + \delta_{m+1}$ then $2^{k_{m+1} -1} \delta_{m+1} -2$
must be divisible by 2.  Therefore it has to be divisible by 6. 
Since it is always divisible by 2 it has to be also divisible by 3.
As before, this shows that the value of $k_{m+1}$ determines the
value of $\delta_{m+1}$ for which this is true.
The corresponding condition takes the form
$$
2^{k_{m+1}} \equiv \delta_{m+1}  \pmod 3 \,.
\tag"{(2.9)}"
$$
From (2.8)
$$
2^{k_{m+1}-1} t-3^{m+1} s_1 = 3r_m^{(1)} 
   - {2^{k_{m+1} -2} \delta_{m+1} -1 \over 3} \,.
$$
A general solution of the last equation is
$t = 3^{m+1} s_2 + \bar{q}_{m+1} , s_1 =2^{k_{m+1} -1} s_2+\bar{r}_{m+1}$
or $s = 2s_1 = 2^{k_{m+1}} s_2 + 2 \bar{r}_{m+1}$.
This gives already one of the statements of the 
Structure Theorem.  For $\bar{r}_{m+1} , \bar{q}_{m+1}$
we have the equation 
$$
2^{k_{m+1}-1} \bar{r}_{m+1} - 3^{m+1} \bar{q}_{m+1}
 = r_m^{(1)} + {1-\delta_{m+1} 2^{k_{m+1} -2}\over 3}
\tag"{(2.10)}"
$$

Now we argue in the same way as in the case of $m=1$. 
If $rhs$ of (2.10) is non--negative, we can always find unique
$\bar{r}_{m+1} , \bar{q}_{m+1} , 1 \leq \bar{r}_{m+1}$
$\leq 3^{m+1} , 1 \leq \bar{q}_{m+1} \leq 2^{k_{m+1}-1}$,
for which (2.10) is true.

If $rhs$ of (2.10) is negative we modify it as follows
$$
2^{k_{m+1} -1} (\bar{r}_{m+1} -1) -
   3^{m+1} \bar{q}_{m+1} = 2^{k_{m+1} -1} +
  r_m^{(1)} + {1-\delta_{m+1} 2^{k_{m+1}-2}\over 3} \,.
$$
Now the $rhs$ is positive and we can find a solution for which 
$1 \leq \bar{q}_{m+1} \leq 2^{k_{m+1}-1}, 1 \leq \bar{r}_{m+1} \leq 3^{m+1}$.
In all cases $ q_m + 2\bar{r}_{m+1} , r_{m+1} = \bar{q}_{m+1}$.

If $r_m$ is odd, $r_m = 2 r_m^{(1)} +1$ and
$\delta_m=1$ then
$$
3 \cdot 3^{m+1} \cdot s+ 18 r_m^{(1)} + 11 = 2^{k_{m+1} -1} y
\tag"{(2.11)}"
$$
and $s$ has to be even, $s= 2s_1 +1$.  This yields
$$
6 \cdot (3^{m+1} s_1 + 3 r_m^{(1)} +2 ) + 3^{m+2} -1 =
   2^{k_{m+1} -1} (6t+\delta_{m+1})
$$
and thus $3^{\2m}-1 -2^{k_{\1m}-1}\delta_{\1m}$ must be
divisible by 6.  Therefore $\delta_{\1m}$ should be such that
$2^{k_{\1m} -1} \delta_{\1m} +1$ is divisible by 3 which
is equivalent to $2^{k_{\1m}} \equiv \delta_{\1m} \pmod 3$.
It is clear that $3^{\2m} -1-2^{k_{\1m}-1} \delta_{\1m}$ is even.

Now we write as before
$$
t = 3^{\1m} s_2 + \bar{q}_{\1m}, s_1 = 
   2^{k_{\1m}-1} s_2 + \bar{r}_{\1m} \,,
$$
and get for $\bar{q}_{\1m} , \bar{r}_{\1m}$ the equation
$$
2^{k_{\1m}-1} \bar{q}_{m+1} -3^{m+1} r_{m+1} =
  3r_m^{(1)} + 2 +
   {3^{\2m}-1 -1^{k_{\1m} -1} \delta_{\1m}\over 6} \,.
$$
This shows that $r_{\1m} = \bar{q}_{\1m} , q_{\1m} = 2 \bar{r}_{\1m} + q_m$.
The case $\delta_m = -1$ is considered in a similar way.
The Structure Theorem is proven.

\head{\bf \S3.  A Corollary of the Structure Theorem}\endhead

Take $x_0 \in \sqcap$ and put $x_m = T^m x_0$,
$y_m = \ell n \; x_m , z_m = y_m - y_0, m \geq 1$.

Consider the probability
$$
P_m (a,b) = P \left\{ x_0 \big| a \leq
{z_m +m (2 \ell n  \; 2 - \ell n \; 3)\over \sqrt{\sigma m}} \leq b \right\} \,.
$$
Here $a, b$ are fixed numbers, 
$\sigma > 0$ is a constant which will be described during the proof,
the probability means the normalized wrt $\sqcap$ density,
i.e. $P_m$ is the relative (wrt $\sqcap$) density of $x_0 \in \sqcap$
satisfying the above inequalities.

\proclaim{Theorem 3.1}
$$
\lim\limits_{m \rightarrow \infty} P_m (a,b)  = {1\over \sqrt{2 \pi}}
   \int\limits_a^b e^{-u^2/2} du \,.
$$
\endproclaim

\demo{Proof}.  Consider any progression
$\Sigma^{(k_1, \ldots, k_m, \epsilon_m)}$
(see the formulation of the Structure Theorem in \S1).
Then its probability in the sense mentioned above
$$
P \{ \Sigma^{(k_1 , \ldots, k_m, \epsilon _m)} \} =
   3\cdot {1\over 6\cdot 2^{k_1 +\cdots+ k_m}} =
  {1\over 2^{k_1 +\cdots + k_m +1}} \,.
\tag"{(3.1)}"
$$
Actually the factor 3 is connected with the normalization
density $(\sqcap$) $=$ density ($\sqcap^{+1}$) $+$ density
($\sqcap^{-1}$) $= {1\over 6} + {1\over 6} = {1\over 3}$
and additional 1 is connected with uniform distribution of
$\epsilon_m = \pm 1$ independent on $k_j$.

Take large $x_0 \in \Sigma^{(k_{1 \cdots} k_m, \epsilon _m)}$,
$x_0= (6( 2^{k_1 +\cdots + k_m}p + q_m ) + \epsilon _m )$.
Then $x_0 = 6p\cdot 2^{k_1 +\cdots +k_m} (1 + \circ (1))$,
$x_m = 6\cdot p \cdot 3^m (1+ \circ (1))$ and
$z_m=\ell n \;{x_m\over x_0}=(k_1+\cdots +k_m )\ell n \;2-m \;\ell n\;3+0(1)$.
Therefore $z_m + m ( 2 \ell n \; 2 - \ell n \; 3) = (k_1 +\cdots + k_m-2 m)$
$\ell n \; 2 (1  + \circ (1))$.
\enddemo
It follows from (3.1) that $k_1,\ldots , k_m$ are independent random 
variables having geometrical distribution with parameter ${1\over 2}$ and
${k_1 + \cdots + k_m - 2m\over \sqrt{\sigma\, m}}$ for some $\sigma >0$
has limiting Gaussian distribution.
This implies the statement of the theorem.  In an analogous way one
can prove a limiting theorem for finite--dimensional distributions of
$z_m$ mentioned in \S1.

Theorem 3.1 says that for very large 
$x_0$ typical $z_m = \ell n {x_m\over x_0}$ 
decrease with the drift coefficient 
$ - (2 \; \ell n\; 2 - \ell n \; 3)$.
This means that $x_m$ also decrease and this gives some reasons to
expect that $(3x +1)$--problem is true.
However, the main difficulty lies in the dynamics on the intermediate
scales.

\head{{\bf\S4}.  $\bold{(r_m},\,\bold{\delta_m )}$ {\bf as a Random Walk}}\endhead

Take $(r_m, \delta_m)$ and a progression 
$\ssize{\Sigma^{(k_1,\ldots ,k_{m}}}, \ssize{\epsilon)}$ such that 
$T^m \Sigma^{(k_1,\ldots ,k_m , \epsilon)} = \wedge^{(r_m, \delta_m)}$ 
(see the Structure Theorem). 
In principle it can happen that for some
$(r_m , \delta_m)$ there are no such 
$\Sigma^{(k_1 ,\ldots ,k_m , \epsilon)}$.
But if they exist then
$T^j \Sigma^{(k_1 , \ldots , k_j , \epsilon )} = \wedge ^{(r_j , \delta_j)}$,
$1 \leq j \leq m$, and the sequence $(r_j , \delta_j)$, 
$1 \leq j \leq m$ can be viewed as a trajectory of some
random walk which ends at $(r_m , \delta_m)$.
Different $\Sigma^{(\kme)}$ generate different trajectories.

We shall use the notation 
$\Phi_m(k_1, k_1, \ldots , k_m , \epsilon)$
$=(r_m (\kme)$, $\delta_m (\kme))$. 
It is clear that
$\Sigma^{(\kme)} \subset \Sigma^{(k_1 , \ldots, k_{m-1}, \epsilon )}$
and $\Phi_j (k_1 , \ldots , k_j , \epsilon) = (\rjdj)$ where
$T^j\Sigma^{(k_1 ,\ldots , k_j , \epsilon)}=\wedge^{(\rjdj)}$,
$1\leq j \leq m$.  Sometimes we shall use also the equivalent
writing $(r_m(\kme)$, $ \delta_m(\kme)) = \Phi_m(q_m(\kme), \epsilon)$.
The value of $\delta_j$ can be found from (2.9):
$$
2^{k_j} \equiv \delta_j \pmod 3
\tag"{(4.1)}"
$$
which imposes some restrictions on possible values of
$k_j$ provided that $\delta_j$ are given.
We shall show that there is another restriction of a similar type.

As in \S2 we have the equation
$$
3[6(3^{m-1} p+ r_{m-1}) +
   \delta_{m-1}]+1 = 2^{k_m}y \, ,
$$
$p\geq 0$ and 
$6(3^{m-1} p+ r_{m-1}) +\delta_{m-1} \in \wedge^{(r_{m-1},\delta_{m-1})}$.

Since $y \in \sqcap$ we write $y = 6s + \delta_m$ where $\delta_m$
is found from (4.1) and
$$
6[2^{k_m} s - 3^m p] + 2^{k_m}\delta_m - 
    3\delta_{m-1} -1 = 18r_{m-1} \,.
$$
Define $t$ by setting $s = 3^mt+r_m$, and then define
$t_m$ by  $p= 2^{k_m}t+t_m$.  Then
$$
6[2^{k_m} r_m - 3^m t_m] = 18 r_{m-1} + 3\delta_{m-1}
   +1-2^{k_m} \delta_m \,.
\tag"{(4.2)}"
$$
(4.2) shows that for given $\delta_{m-1}$ the value of
$k_m$ should be such that
$$
3 \delta_{m-1} + 1 \equiv 2^{k_m} \cdot \delta_m \quad \pmod 6
$$
or $2^{k_m -1} \cdot \delta_m \equiv {3\delta_{m-1} +1\over 2} \pmod 3$.
Using (4.1) we can write
$$
2^{k_m} = \delta_m + 3 a_m^{(1)}
$$
for some odd $a_m^{(1)}$.  Then
$$
{2^{k_m} \delta_m - 3 \delta_{m-1} -1\over 6} =
   {a_m^{(1)} \delta_m - \delta_{m-1}\over 2} \;.
$$
Returning back to (4.2) we get
$$
2^{k_m} r_m  - 3^m t_m = 3r_{m-1} -
   {a_m^{(1)} \delta_m - \delta_{m-1}\over 2 }\; .
\tag"{(4.3)}"
$$
This shows that for given $r_m$ the value of $k_m$ should
be such that
$$
2^{k_m} r_m + {a_m^{(1)} \delta_m - \delta_{m-1}\over 2}
    \equiv 0 \quad \pmod 3 \; .
\tag"{(4.4)}"
$$
Since $a_m^{(1)}$ is odd, $a_m^{(1)} = 2a_m^{(2)} +1$
and $a_m^{(2)} = g_m + 3a_m^{(3)}$.  Remark that 
$a_m^{(1)}$, $a_m^{(2)}$, $a_m^{(3)}$, $g_m$ are
functions of $k_m$ only.  Actually
$$
2^{k_m} = \delta_m + 3 + 6 g_m + 18 a_m^{(3)} \;.
\tag"{(4.5)}"
$$
For $r_m$ we can write
$r_m = h_m +3r_m^{(1)}$ where $h_m$ can take values
$0, 1, 2$.  The last expression can be considered as
the definition of $h_m , r_m^{(1)}$ as functions of
$r_m$.  From (4.4)
$$
h_m \delta_m + g_m \delta_m +
{\delta_m - \delta_{m-1}\over 2} \equiv 0 \quad \pmod 3
\tag"{(4.6$^\prime$)}"
$$
or
$$
h_m + g_m + {1-\delta_{m-1} \delta_m\over 2}
    \equiv 0 \quad \pmod 3 \;.
\tag"{(4.6$^{\prime\prime}$)}"
$$

The equations (4.6$^\prime$), (4.6$^{\prime\prime} $) have an
important interpretation.  Suppose that we are given
$r_m$, $\delta_m$, $\delta_{m-1}$.  Then the value of $\delta _m$
determines the parity of $k_m$, the value of $r_m$ gives the
value of $h_m$ and (4.6$^{\prime\prime}$) allows us to find the
value of $g_m$.

Take again (4.3).  It shows how to find $r_{m-1}$ knowing 
$r_m$, $k_m$, $\delta_{m-1}$.
From (4.5), (4.6$^\prime$), (4.6$^{\prime\prime}$)
$$
2^{k_m} r_m^{(1)} + a_m^{(2)} \delta_m +
  {\delta_m - \delta_{m-1}\over 2}
    -3^{m-1} t_m = r_{m-1} \;.
\tag"{(4.7)}"
$$
Using the analogy with Markov processes we can call
(4.7) the backward system of equations.

\head{\bf \S5.  The Ensemble} $\bold{\Phi_m^{-1}}$$\bold{(r_m ,\delta_m )}$.\endhead

As it follows from \S3 it is natural to consider
the probability distribution $P$ for which $\epsilon = \pm 1$
with probabilities $\frac{1}{2}$ and
$k_1 , k_2 , \ldots , k_m , \ldots$ is a sequence
of independent, random variables, also independent on
$\epsilon$ and having the geometric distribution with
exponent $\frac{1}{2}$.  All probabilities which we
consider below are induced by this distribution.
For example, with respect to this distribution the
probability of any $q_m (k_1, \ldots , k_m, \epsilon)$
equals $\frac{1}{2^{k_1 + \cdots + k_m +1}}$
and the probability of a pair $(r_m , \delta_m)$ is the 
probability of all $(q_1 , \ldots , q_m, \epsilon)$
which give $(r_m, \delta_m)$ under the mapping $\Phi_{m}$.
The main purpose of this section is to study the
probabilities of pairs $(r_m , \delta_m)$.

The pair $(r_m, \delta_m)$ can take 2.3$^m$ values.
On the other hand the number of typical $(\kme )$ grows
(in a weak sense) as $2^{2m}$.  Therefore it is 
natural to expect that typically $\Phi_m^{-1} (r_m, \delta_m)$
contains $2^{2m} \cdot 3^{-m}$ elements. 

Put $-c(k_m, \delta_m , \delta_{m-1}) = {a_m^{(1)} \delta_m - \delta_{m-1}\over 2}$, where (see above) $a_m^{(1)} = {2^{k_m} - \delta _m\over 3}$.
Thus
$-c(k_m, \delta_m, \delta_{m-1}) = {2^{k_m} \delta_m - 1 - 3\delta_{m-1}\over 5}$.
In particular,
$-c (1, -1, \delta_{m-1}) = - {1+\delta_{m-1}\over 2}$,
$-c(2, 1, \delta_{m-1}) = {1-\delta_{m-1}\over 2}$,
and so on.  It is clear that
$c(k_m, \delta_m , \delta_{m-1})$ can be positive or negative
and for large $k$
$$
-c(k_m, \delta_m , \delta_{m-1}) \sim 
   {2^{k_m} \delta_m\over 6} \,.
$$
Denote $\rho_m = {r_m\over 3^m}$.  Then $0 \leq \rho_m \leq 1$
and possible values of $\rho_m$ go with the step ${1\over 3^m}$.
From (4.3)
$$
\rho_m = {t_m\over 2^{k_m}} + {1\over 2^{k_m}} \cdot
    \rho_{m-1} + {c(k_m , \delta_m , \delta_{m-1})\over 2^{k_m} \cdot 3^m} \;.
\tag5.1
$$
The iteration of the last equality yields
$$
\rho_m = {t_m\over 2^{k_m}} + {t_{m-1}\over 2^{k_m + k_{m-1}}}
+\cdots + {t_1\over 2^{k_m+k_{m-1} +\cdots +k_1}} +
  \sum\limits_{s=1}^m
     {c(k_s , \delta_s , \delta_{s-1})\over
     2^{k_m +\cdots + k_s} \cdot 3^s} \;.
\tag5.2
$$
It follows easily from (4.3) and from \S2 that if
$q_m = q_m (k_1 , \ldots , k_m , \epsilon) \in \Phi_m^{-1} (r_m, \delta_m)$
then
$$
q_m = t_m \cdot 2^{k_{m-1}+\cdots + k_1} +
   t_{m-1} 2^{k_{m-2}+\cdots +k_1} + \cdots + 
     t_2 \cdot 2^{k_1} + t_1 \;. 
\tag5.3
$$
Put $\kappa_m = q_m 2^{-(k_m +\cdots + k_1)}$.  We have
$$
\kappa_m = q_m 2^{-(k_m + \cdots + k_1)} =
  {t_m\over 2^{k_m}} + {t_{m-1}\over 2^{k_m + k_{m-1}}}
    +\cdots + {t_1\over 2^{k_m +k_{m-1} + \cdots + k_1}}
\tag5.4
$$
and from (5.2)
$$
\rho_m = \kappa_m + \sum\limits_{s=1}^m
   {c(k_s, \delta_s , \delta_{s-1})\over 2^{k_m + \cdots + k_s} \cdot 3^s} =
  \kappa_m + {1\over 3^m} \sum\limits_{s=1}^m
     {3^{m-s} c(k_s , \delta_s , \delta_{s-1})\over 2^{k_m +\cdots+ k_s}} \;.
\tag5.5
$$
Since $k_s$ are independent random variables having
geometric distribution with parameter ${1\over 2}$,
$\delta_s = -1$ or $+1$ depending on the parity of $k_s$
the sum $k_m +\cdots + k_s$ grows typically as
$2(m-s)$.  By this reason the last sum in (5.5) is converging,
at least in probability,
takes values $O(1)$ and has limiting distribution as 
$m \rightarrow \infty$.  The formula (5.5) shows that
for $(k_1 , \ldots , k_m , \epsilon ) \in \Phi_m^{-1} (r_m , \delta_m)$
the difference $\rho_m - \kappa_m = O({1\over 3^m})$.
Write $k = k_1 + \cdots + k_m$.
It is a well--known combinatorial fact that the number
$H_m(k)$
of solutions of the last equation with $k_i \geq 1$ equals to
$$
H_m (k) = 
\pmatrix k-1 \\ m-1\endpmatrix = 2^{k-1} \cdot
 G_m(k-2m)
\tag5.6
$$
where $G_m(k-2m)$ have Gaussian asymptotics
$$
G_m (k - 2m) \sim {1\over \sqrt{2\pi \sigma m}}
   \exp \biggl\{-{(k-2m)^2\over 2 \sigma m} \biggr\}
$$
for some constant $\sigma > 0$ and not too large
$|k-2m|$.

Put
$$
\theta_m = \sum\limits_{s=1}^m 
  {3^{m-s} c(k_s , \delta_s , \delta_{s-1})\over
     2^{k_m +\cdots + k_s}}
$$
and
$$
A_{m,i} = \bigg\{\big((k_1, \ldots , k_m ), \epsilon \big) \bigg|
  {i\over 10} \leq \theta_m < {i+1\over 10} \bigg\} \;.
$$
Instead of 10 we could take any large enough integer.
It is clear that the value of $\theta_m$ is basically determined by
the last $k_m, k_{m-1} , \ldots \; .$
It follows from (5.5) that 
$((k_1, \ldots , k_m), \epsilon) \in A_{m,i}\cap \Phi_m^{-1}((r_m, \delta_m))$
iff
$$
\rho_m - {(i+1)\over 10 \cdot 3^m} <
   \kappa_m \leq \rho_m - {i\over 10 \cdot 3^m} \;.
\tag5.7
$$
It is easy to show that one can find such constant
$\gamma_\circ > 0$ that
$$
P\{|\theta_m| > m^{\gamma_\circ} \} \leq {1\over m} \;.
$$
We shall use the notation $D_m^\prime$ for the set of
$(k_1 , \ldots k_m , \epsilon )$ for which $|\theta_m| \leq m ^{\gamma_\circ}$.

For any value of $k$ the number of possible 
$(k_1 ,\ldots , k_m , \epsilon ) \in \Phi_m^{-1}(r_m ,\delta_m )\cap A_{m,i}$
with the given $k$ is at most  ${2^k\over 10.3^m}$
because the interval (5.7) has width $1/10 \cdot 3^m$ and each 
$\kappa_m$ is rational with denominator $2^k$ and all $\kappa_m$ are distinct
by (5.3).

Therefore the probability
of this set is not greater than 
${1\over 2\cdot 10\cdot 3^m} = {1\over 20\cdot 3^m}$.
As was mentioned above 
$P\{(r_m, \delta_m) \} = \sum_{\Phi_m (q_m , \epsilon) = (r_m , \delta_m)}$
$P \{(q_m , \epsilon)\}$.
Actually we can consider the partition $\xi_m$
of the space $\Omega_m$ of pairs $(\kappa_m, \epsilon )$ onto
pre--images $\Phi_m^{-1} ((r_m , \delta_m))$.
Denote by $H_m$ the entropy of this partition, i.e. 
$H_m = -\sum P ((r_m , \delta_m)) \ln P((r_m, \delta_m))$.
Below the letter $H$ is used for the entropy of a partition.

\proclaim{Theorem 5.1}
$H_m \geq m \ln 3 -(2\gamma_\circ +7) \ln m$
\endproclaim

\demo{Proof}
The proof is based upon the fact that if the entropy is small then there
should be elements of partition having a big measure.  This is
impossible in our case.
Let $B_k = \{(k_1 , \ldots , k_m, \epsilon ) | k_1 + \cdots + k_m = k\}$.
It follows easily from the combinatorial formula above that we can find a 
constant $\gamma_1$ for which for all sufficiently large $m$

$$
P  \biggl\{\ \underset\dsize{ |k-2m| \geq \gamma_1 \sqrt{m \ln m}}\to
   {\hskip -1.2cm \cup \; B_k} \;\;\biggr\}
\leq {1\over m^{\gamma_\circ +2}} \;.
$$
Introduce the partition $\alpha_m$ which has two elements
$$ 
C_m^\prime =  \underset\dsize{ |k-2m| \leq \gamma_1 \sqrt{m \ln m}}\to {\hskip -1.2cm \cup \; B_k} \;,
 \qquad   C_m^{\prime\prime} =  \underset\dsize{ |k-2m| > \gamma_1 \sqrt{m \ln m}}\to {\hskip -1.2cm \cup \; B_k} \;,
$$
and another partition $\beta_m$ also onto two elements,
$D_m^\prime$,$ D_m^{\prime\prime}$ where $D_m^{\prime\prime}$
 is the complement of $D_m^\prime$.  Then
$$
\align
H_m \geq H (\xi_m |\alpha_m \vee \beta_m )& =
-\biggl(\sum\limits_{(r_m, \delta_m)} P ((r_m, \delta_m) \big|
  C_m^\prime \cap D_m^\prime)\biggr. \\
    &\biggl.\qquad \ln P ((r_ m, \delta_m )| C_m^\prime \cap D_m^\prime)
       \cdot P (C_m^\prime \cap D_m^\prime )+ \ldots \biggr)
\tag5.8
\endalign
$$
where dots mean similar sums multiplied by small probabilities
$P (C_m^\prime \cap D_m^{\prime\prime} )$, 
$P (C_m^{\prime\prime}\cap D_m^\prime )$,
$P (C_m^{\prime\prime} \cap D_m^{\prime\prime})$
respectively.  All conditional entropies are less than
$m \ln 3 + \ln 2$ because the partition has not more than
$2.3^m$ elements.  Therefore because of the estimates
of the measures all these terms in (5.7) have absolute values
less than a constant.  Assume that the first sum is smaller than
$m \ln 3 - (2 \gamma_\circ +6) \ln m$.
By Chebyshev inequality
$$
\align 
&P \bigl\{ - \ln P (\Phi_m^{-1} (r_m, \delta_m)
 | C_m^\prime \cap D_m^\prime) \geq m \ln 3 -
 (2 \gamma_\circ +3) \ln m \bigr\} \\
& \qquad \leq {m \ln 3 -(2 \gamma_\circ +6) \ln m \over m \ln 3 
-(2\gamma_\circ +3) \ln m}=1- 
     {(2 \gamma_\circ +3) \ln m\over m \ln 3 - (2\gamma_\circ +3) \ln m} \;.
\endalign
$$
Therefore
$$
\align
P \bigl\{ - \ln P
  \big(\Phi_m^{-1} (r_m , \delta_m)\big|
    C_m^\prime \cap D_m^\prime \big) 
     & < m \ln 3 - (2 \gamma_\circ +3 ) \ln m \bigr\} \\
& \geq {(2 \gamma_\circ +3) \ln m\over
    m \ln 3 - (2\gamma_\circ +3) \ln m}
\endalign
$$
and by this reason the set of $(r_m , \delta_m)$ for which 
$-\ln P (\Phi_m^{-1}(r_m ,\delta_m)|C_m^\prime \cap D_m^\prime)$
$<m \ln 3 - (2 \gamma_\circ +3) \ln m$ or, equivalently, 
$P (\Phi_m^{-1} (r_m , \delta_m) | C_m^\prime \cap D_m^\prime)$
$> {m^{2 \gamma_\circ +3}\over 3^m}$ is not empty. 
We shall show that
this is impossible. 

By definition
$$
\align
P\bigl(\Phi_m^{-1} (r_m , \delta_m) \big| C_m^\prime \cap D_m^\prime
\bigr\} &= {P \bigl\{\Phi_m^{-1} (r_m , \delta_m) \cap
    C_m^\prime \cap D_m^\prime \bigr) \over
      P (C_m^\prime \cap D_m^\prime )}\\
&= {1\over P(C_m^{\prime} \cap D_m^\prime )}
  \sum\Sb |k-2m|\leq \gamma_1 \sqrt{m \ln m}\\
          |i| \leq m^{\gamma_\circ} \endSb \\
& \qquad P \bigl( \Phi_m^{-1} (r_m, \delta_m) \cap A_{m,i}
     \cap B_k \bigr) \;.
\endalign
$$
Therefore one can find $i_\circ$, $k_\circ$ such that
$$
P\bigl(\Phi_m^{-1} ( r_m , \delta_m ) \cap
     A_{m,i_0} \cap B_{k_0} \bigr) \geq {m^2\over 3^m}
$$
for all sufficiently large $m$.
But it was already shown above that this probability cannot be
greater than ${1\over 20.3^m}$.  This contradiction proves
the theorem.
\enddemo

Theorem 5.1 shows in what sense the distribution
$\bigl\{P \bigl(\Phi_m^{-1} (r_m, \delta_m)\bigr)\bigr\}$
is close to the uniform.
We believe that actually $H_m \geq m \,\ln\; 3-$const.


\Refs

\widestnumber\key{999}


\ref \key{BP} \by  K. A. Borovkov and D. Pfeifer    \pages 300--310
\paper Estimates for the Syracuse Problem via a probabilistic model
\yr 2000 \vol 45, N2
\jour Theory of Probability and its Applications \endref

\ref \key{L} \by J. C. Lagarias \pages 3--23
\paper The $(3x+1)$--problem and its generalizations
\yr 1985 \vol 92
\jour American Mathematical Monthly \endref


\ref \key{W} \by G. J. Wirsching
\paper  The Dynamical System generated by the $(3x+1)$--function
\inbook Lecture Notes in Mathematics, N1681
\publ Springer--Verlag
\publaddr Berlin \yr 1998,  158p \endref
\end{document}